\documentclass[a4paper]{amsart} 
\usepackage{hyperref}

\theoremstyle{plain} 
 \newtheorem{thm}{Theorem} 
 \newtheorem{cor}[thm]{Corollary}
 
 \newtheorem{lem}[thm]{Lemma}
\theoremstyle{definition} 
 
\theoremstyle{remark} 
\newtheorem*{rem}{Remark}

\newcommand{\we}{\wedge} 
\newcommand{\der}{\Delta} 
\renewcommand{\d}{d} 
\newcommand{\del}{\partial} 
\newcommand{\dlog}[1]{\dfrac{\d{#1}}{#1}} 
\newcommand{\dlogg}[1]{{\d{#1}/{#1}}} 

\newcommand{\A}{{\mathcal A}} 
\newcommand{\C}{{\mathcal C}} 
\DeclareMathOperator{\sign}{sign} 
\newcommand{\dc}{\operatorname{\mathbf{d}}} 
\newcommand{\den}{\mathrm{D}} 

\begin{document} 

\title[Non-vanishing of the twisted cohomology]{ 
Non-vanishing of the twisted cohomology \\ 
on the complement of hypersurfaces } 
\author{Yukihito Kawahara} 
\address{Department of Mathematics, Tokyo Metropolitan University \\ 
Minami-Ohsawa 1-1, Hachioji-shi, Tokyo  192-0397, Japan } 
\email{ykawa@comp.metro-u.ac.jp} 
\maketitle 
\begin{abstract} 
Under the generic situation, 
the cohomology with the coefficients in the local system  
on complements of hypersurfaces 
vanishes except in the highest dimension. 
Our problem is of when the local system cohomology does not vanish. 
In the case of arrangements of hyperplanes, many examples were founded.   
In this paper, we shall generalize their examples to hypersurfaces. 
We obtain that  hypersurfaces given by some linear system 
have non-vanishing local system cohomologies. 
\end{abstract} 

\section{Introduction} \label{sec1} 

Let $V_1, \ldots, V_m$ be hypersurfaces 
in the complex projective space $\mathbb{P}^n$ of dimension $n$ 
and let $\mathcal{L}$ be a complex local system of rank one over 
the complement $M = \mathbb{P}^n \setminus \cup_{i=1}^{m} V_i$. 
In \cite{A, KN, Ch, Di}, we know vanishing theorems of 
the twisted cohomology that is  
the cohomology with the coefficients 
in the local system $\mathcal{L}$ over $M$  as follows.   
Under the generic situation for $V_1, \ldots, V_m$, 
if $\mathcal{L}$ is non-trivial and generic,  
then the twisted cohomology on $M$ 
vanishes except in the highest dimension: 
$$   H^k(M, \mathcal{L}) = 0  \qquad \text{ for } k \not= n.   $$ 
In particular, 
vanishing theorems for the case of hyperplanes was found 
in \cite{Ko, Yu, CDO, Ka2} (cf. \cite{ESV, STV}). 
Recently, 
arrangements of hyperplanes with non-vanishing twisted cohomology  
were studied well and 
many examples were found (cf. \cite{CS, Fa, LY, Ka3}). 
\cite{Yu2} implied that 
most of them  consist of special elements of pencils in $\mathbb{P}^2$. 
In this paper 
we generarize it to hypersurefaces in the general dimension. 
We shall show that 
hypersurfaces that is the support of some divisors in some linear system 
have non-vanishing twisted cohomologies besides the top dimension.  

Let $\Omega_M$ denote the sheaf of germs of holomorphic forms on $M$ 
and let $\mathcal{O}_M = \Omega_M^0$. 
Let $D_1, \ldots, D_m$ be effective divisors on $\mathbb{P}^n$  
such that the support of $D_i$ is $V_i$ for $i$ and 
$D_i$ and $D_j$ are linearly equivalent for $i \not= j$. 
In this case, $D_i$'s have the same degree (see \cite{H}).  
Let $\lambda = (\lambda_1, \ldots, \lambda_m)$ be 
a complex weight with $\sum_{i=1}^{m} \lambda_i =0$. 
For $i \not= j$, 
we have $D_i - D_j$ is the divisor of some rational function $f_{ij}$. 
Fix $j$ and define the global one-form 
$\omega_{\lambda} = \sum_{i \not= j} \lambda_i \d{\log{f_{ij}}}$.    
This is independent on the choice of $j$ 
and then we can denote       
$$   \omega_{\lambda} = \sum_{i=1}^{m} \lambda_i \d{\log{D_i}}   
        \in \Gamma{(M, \Omega_M^1)},  
  \qquad  \sum_{i=1}^{m} \lambda_i =0.   $$  
We can define the flat connection 
$\nabla_{\lambda} = d + \omega_{\lambda} \we : \mathcal{O}_M \to \Omega_M^1$ 
and define the local system $\mathcal{L}_{\lambda}$ by its kernel. 
Since $M$ is a Stein manifold, we have 
$$  H^k(M, \mathcal{L}_{\lambda}) 
     \simeq H^k( \Gamma{(M, \Omega_M)}, \nabla_{\lambda})  $$  
and  $H^k(M, \mathcal{L}_{\lambda}) =0$ for $k >n$ (cf. \cite{De}).

\begin{thm}\label{MainTheorem} 
Let $1 < n < s \leq m$. 
Let $D_1, \ldots, D_m$ be effective divisors on $\mathbb{P}^n$ 
with same degree  
and let $M$ be the complement of their supports.  
Suppose 
\begin{itemize}
	\item[(A1)] 
$D_1, \ldots, D_s$ are elements 
of some linear system $\Lambda$ on $\mathbb{P}^n$ 
with dimension $n-1$.  
	\item[(A2)] 
There exist a base point $P$ of $\Lambda$ 
such that $D_1 \cdots D_n$ is normal crossing locally at $P$ 
and $D_i$ does not pass through $P$ for $i = s+1, \ldots, m$. 
	\item[(A3)] 
$D_1, \ldots, D_s$ are 
in general position as points in $\Lambda = \mathbb{P}^{n-1}$.  
\end{itemize} 
If $\lambda$ is a non-trivial weight, that is, $\lambda \not\in \mathbb{Z}^m$, 
such that $\sum_{i=1}^s \lambda_i =0$ and 
$\lambda_i =0$ for $i = s+1, \ldots, m$, 
then we have 
$$   \dim{H^{n-1}(M, \mathcal{L}_{\lambda})} \geq  \binom{s-2}{n-1}.   $$ 
Moreover, if $s < m$ then we have 
$$   \dim{H^{n}(M, \mathcal{L}_{\lambda})} \geq  \binom{s-2}{n-1}.   $$ 
\end{thm}  
\begin{rem}  
If the support of $\sum_{i=1}^m D_i$ contains a hyperplane $H$, 
we can consider $M$ as the complement of affine hypersurfaces 
in $\mathbb{C}^n = \mathbb{P}^n \setminus H$. 
The setting of the affine hypersurfaces case was found in \cite{KN, Ki} 
(cf. Section \ref{affine}).  
\end{rem} 
\begin{rem}  
Let $D_1, \ldots, D_m$ be 
divisors of degree $d$ defined only by hyperplanes. 
Let $\A$ be the set of irreducible components of the support 
of $\sum_{i=1}^m D_i$. 
Then $\A$ is an arrangement of hyperplanes in $\mathbb{P}^n$ and 
$M= \mathbb{P}^n \setminus \cup_{H \in \A} H$. 
If $D_1, \ldots, D_m$ are divisors of degree one, then 
they are defined by hyperplanes and 
this result is known (cf. \cite{Fa}). 
\end{rem}

\section{Preliminary: Rational forms} \label{preliminary} 

Let $V$ be a vector space of dimension $\ell$  
over a field $K$ of characteristic zero.  
Let $S = K[V]$ be the symmetric algebra of $V^*$ 
and let $F = K(V)$ be the quotient field of $S$.   
We consider $S$ as the polynomial algebra 
and $F$ as  the field of rational functions on $V$. 
Let $\Omega{(V)} = \oplus_{p=0}^{\ell} \Omega^p{(V)}$ be  
the exterior algebra of $F$-vector space $F \otimes V^*$ and 
let $d$ be the usual differential.  
When we choose a basis $x_1, \ldots, x_{\ell}$ of $V^*$, 
we have 
$S = K[x_1, \ldots, x_{\ell}]$,  $F = K(x_1, \ldots, x_{\ell})$, 
$$  \d{f} = \sum_{i=1}^{\ell} \dfrac{\del{f}}{\del{x_i}} \otimes x_i 
          = \sum_{i=1}^{\ell} \dfrac{\del{f}}{\del{x_i}} \d{x_i}  
     \quad \text{ for }  f \in F,  $$  
$\Omega^0{(V)} = F$,  
$\Omega^1{(V)} = F \otimes V^* 
   = F \d{x_1} \oplus \cdots \oplus F \d{x_{\ell}}$,  and, 
$\Omega^p{(V)} = \oplus_{i_1 < \cdots < i_p} 
     F \d{x_{i_1}} \we \cdots \we \d{x_{i_p}}$.  
For $p \geq 2$ and  $\omega_1, \ldots, \omega_p \in \Omega^1{(V)}$,  define 
$$  \der{\left[ \omega_1 : \cdots : \omega_p \right]} := 
 \sum_{k=1}^{p} (-1)^{k-1} \omega_1 \we \cdots \we \widehat{\omega_k} \we 
        \cdots \we \omega_p.   $$ 

\begin{lem}\label{lemma1}  
Let $p \geq 2$ and  $\omega_1, \ldots, \omega_p \in \Omega^1{(V)}$. 
\begin{enumerate}
	\item 
For a permutation $\sigma$ of $\{ 1, \ldots, p \}$, 
we have 
$$  \der{\left[ \omega_{\sigma(1)} : \cdots : \omega_{\sigma(p)} \right]} 
   = \sign{(\sigma)} 
 \der{\left[ \omega_1 : \cdots : \omega_p \right]}.  $$  
	\item 
If $2 \leq j \leq  p-2$, then 
\begin{align*}
  \der{[ \omega_1 : \cdots : \omega_p ]} =  &   
      \der{[ \omega_1 : \cdots : \omega_j ]} 
         \we \omega_{j+1} \we \cdots \we \omega_p  \\ 
    &  + (-1)^j \omega_1 \we \cdots \we \omega_j \we 
     \der{[\omega_{j+1} : \cdots : \omega_p ]}.  
\end{align*}
	\item 
$\der{\left[ \omega_1 : \cdots : \omega_p \right]} = 
    -  (\omega_1 - \omega_2) \we  
        \der{\left[ \omega_2 : \cdots : \omega_p \right]}$.  
	\item 
$\der{\left[ \omega_1 : \cdots : \omega_p \right]} = 
    (-1)^{p-1} 
    (\omega_1 - \omega_2) \we (\omega_2 - \omega_3) \we \cdots 
    \we (\omega_{p-1} - \omega_p)$.    
	\item 
$\der{\left[ \omega_1 : \cdots : \omega_p \right]} = 
    (\omega_2 - \omega_1) \we (\omega_3 - \omega_1) \we \cdots 
    \we (\omega_p - \omega_1)$.    
	\item 
$ \omega_1 \we \der{\left[ \omega_1 : \cdots : \omega_p \right]} = 
    \omega_1 \we \omega_2 \we \cdots  \we \omega_p$. 
\end{enumerate}
\end{lem} 

Since (6), if $\der{\left[ \omega_1 : \cdots : \omega_p \right]} =0$ 
then $\omega_1, \ldots, \omega_p$ is $F$-linearly dependent.  
However, the inverse is not true in general.   
The rational function  
$f \in F$ is said to be \textit{homogeneous of degree $d$}, 
if $f= g_1/g_2$ with homogeneous polynomials $g_1$, $g_2$  
and $d= \deg{g_1} - \deg{g_2}$. 
\begin{lem} \label{lem:loghom} 
Assume $p \geq 2$ and $f_1, \ldots, f_p \in F \setminus K$ 
are homogeneous of same degree.  
Then 
$\der{\left[ \dlogg{f_1} : \cdots : \dlogg{f_p} \right]} =0$,  
if and only if, $\dlogg{f_1} \we \cdots \we \dlogg{f_p} =0$. 
\end{lem} 
\begin{proof} 
Let $\omega_1, \ldots, \omega_p \in \Omega^1{(V)}$. 
By the direct computation, we can obtain 
\begin{itemize}
	\item[(1)] 
$\der{\left[ \omega_1 : \cdots : \omega_p \right]} =0$, if and only if, 
there exists 
$(g_1, \ldots, g_p) \in F^p \setminus \{ (0, \ldots, 0) \} $ 
such that $g_1 + \cdots + g_p =0$ and 
$g_1 \omega_1 + \cdots + g_p \omega_p =0$. 
\end{itemize}
Take $\omega_i = \dlogg{f_i}$ and $g_i = h_i f_i$.  
Then we get 
\begin{itemize}
	\item[(2)] 
For $p \geq 2$ and $f_1, \ldots, f_p \in F \setminus K$, we have 
$\der{\left[ \dlogg{f_1} : \cdots : \dlogg{f_p} \right]} =0$, 
if and only if, 
there exists 
$(h_1, \ldots, h_p) \in F^p \setminus \{ (0, \ldots, 0) \}$ 
such that $h_1 f_1 + \cdots + h_p f_p =0$ and 
$h_1 \d{f_1} + \cdots + h_p \d{f_p} =0$. 
\end{itemize}
When $f_i$'s are homogeneous of same degree,  by using the Euler derivation, 
$h_1 \d{f_1} + \cdots + h_p \d{f_p} =0$ induces 
$h_1 f_1 + \cdots + h_p f_p =0$.   
\end{proof} 
In particular, 
if $c_1 f_1 + \cdots + c_p f_p =0$  
for some $(c_1, \ldots, c_p) \in K^p \setminus \{ (0, \ldots, 0) \}$, 
then 
$\der{\left[ \dlogg{f_1} : \cdots : \dlogg{f_p} \right]} =0$. 
If  $f_1^{n_1} \cdots f_p^{n_p} =1$ 
for some non-zero integers $n_1, \ldots, n_p$ 
then 
$\der{\left[ \dlogg{f_1} : \cdots : \dlogg{f_p} \right]} =0$ also. 
\begin{rem} 
$\der$ is a natural generalization of the linear derivation 
on Orlik-Solomon Algebras (\cite{OT}), which is in the degree one case. 
\end{rem}

\section{Proofs} \label{proofs} 

\begin{proof}[Proof of Theorem \ref{MainTheorem}]  
Let $[x_0: x_1: \ldots: x_n]$ be homogeneous coordinates of $\mathbb{P}^n$.  
We can assume  
$D_i$ is given by a homogeneous polynomial $F_i(x)$ of degree $d$.  
So we can write 
$\omega_{\lambda} = \sum_{i=1}^{m} \lambda_i \dlogg{F_i}$. 
It is easy to check that 
$\dlogg{F_i} - \dlogg{F_j}$ and $\omega_{\lambda}$ are global forms. 
For $1 \leq i_1, \ldots, i_p \leq m$, 
define a holomorphic form on $M$ by  
$$  \eta{[i_1, \ldots, i_p]} 
     :=  \der{\left[ \dlog{F_{i_1}} : \cdots : \dlog{F_{i_p}}  \right]}. $$   
Since Lemma \ref{lemma1} (4), 
if it is not zero then   
$\eta{[i_1, \ldots, i_p]}$ is a global ($p-1$)-form. 
By (A1) and (A3), 
$F_1, \ldots, F_n$ becomes a basis of 
the vector subspace of $\Gamma{(\mathbb{P}^n, \mathcal{O}(d))}$ 
defining the linear system $\Lambda$. 
So we can write 
$F_j = a_{1j} F_1 + \cdots + a_{nj} F_n$ for some constant $a_{ij}$ and 
define the $n \times s$-matrix $A = (a_{ij})$. 
By (A3),  any $n \times n$-minor of $A$ is not zero. 
Due to Lemma \ref{lem:loghom}, we have $\eta{[i_1, \ldots, i_{n+1}]} =0$ 
for $1 \leq i_1 < \cdots < i_{n+1}  \leq s$.  
Because we can write 
$\omega_{\lambda} = \sum_{i \not= i_1} \lambda_i 
  ( \dlogg{F_i}- \dlogg{F_{i_1}})$, 
using Lemma \ref{lem:loghom},  we have 
$\omega_{\lambda} \we \eta{[i_1, \ldots, i_{n}]} =0$ 
and then 
$\nabla_{\lambda}(\eta{[i_1, \ldots, i_{n}]}) =0$ 
for $1 \leq i_1 < \cdots < i_n  \leq s$. 
We note that, 
for $1 \leq i_1 < \cdots < i_{n} \leq s$,   
by (A2) and (A3),  we have 
$\dlogg{F_{i_1}} \we \cdots \we \dlogg{F_{i_n}} \not=0$ 
and, by Lemma \ref{lem:loghom}, 
$\eta{[i_1, \ldots, i_n]} \not= 0$.   
Thus, 
$\eta{[i_1, \ldots, i_{n}]}$ is a $\nabla_{\lambda}$-closed ($n-1$)-form 
for $1 \leq i_1 < \cdots < i_{n} \leq s$. 

By (A2), 
we take a local neighborhood $U$ and coordinates $x=(x_1, \ldots, x_n)$ at $P$ 
such that $D_i$ is defined by $x_i =0$ for $i= 1, \ldots, n$.  
Let $\alpha_j(x) = a_{1j} x_1 + \cdots + a_{nj} x_n$ 
and $H_j = \{ \alpha_j(x) =0 \}$.  
So we get the central arrangement 
$\A = \{ H_j \}_{1 \leq j \leq s}$  
of hyperplanes in $\mathbb{C}^n \simeq U$, 
($\cap_{i=1}^s H_i$ is the origin). 
Let $M(\A)$ denote the complement of an arrangement $\A$. 
Then we have  
$H^k(U \cap M, \mathcal{L}_{\lambda}\vert_{U \cap M}) = 
H^k(M(\A), \tilde{\mathcal{L}}_{\lambda})$ 
where $\tilde{\mathcal{L}}_{\lambda}$ 
is the rank one local system on  $M(\A)$ 
whose monodromy around the hyperplane $H_j$ is 
$\exp{(-2 \pi \sqrt{-1} \lambda_j)}$.  
Since $\lambda$ is non-trivial and $\sum_{i=1}^s \lambda_i=0$, 
without loss of generality, 
we may assume that $\lambda_1$ and $\lambda_s$ are not integers. 
Now, choosing $H_1 \in \A$,  
we get the deconing $\dc{\A}$ (see \cite{OT}), which is 
an arrangement of $s-1$ affine hyperplanes in $\mathbb{C}^{n-1} \simeq H_1$.   
Note that $M(\A) \simeq M(\dc{\A}) \times \mathbb{C}^*$ 
by the restriction of the Hopf bundle.  
Since any $n \times n$-minor of $A$ is not zero, $\A$ is generic and 
$\dc{\A}$ is in general position (\cite{OT}). 
On the other hand, 
for the integer weight $k \in \mathbb{Z}^n$ with $\sum_{i=1}^s k_i =0$, 
we know that   
the local system $\tilde{\mathcal{L}}_{\lambda}$  is equivalent to 
the local system $\tilde{\mathcal{L}}_{\lambda+k}$ associated to 
the integer shift weight $\lambda +k$ (see \cite{OT2}). 
By shifting a weight if necessary,  
we can assume that $\lambda \notin (\mathbb{Z} \setminus \{0\})^n$.    
Since 
$H^k(M(\A), \tilde{\mathcal{L}}_{\lambda}) \simeq 
 H^k(M(\dc{\A}), \tilde{\mathcal{L}}_{\lambda}) \oplus  
 H^{k-1}(M(\dc{\A}), \tilde{\mathcal{L}}_{\lambda})$ (cf. \cite{Fa}),      
in this case, the following is known.  
\begin{lem}[cf. \cite{Ha,KN,Ki,Ka}]  
Let $\A = \{ H_j = \{ \alpha_j =0 \} : 1 \leq j \leq s  \}$ 
be a generic arrangement of hyperplanes in $\mathbb{C}^n$ 
and let $M(\A)$ be its complement. 
For a complex weight 
$\lambda = (\lambda_1, \ldots, \lambda_s)$  
such that $\lambda_1 \notin \mathbb{Z}$, $\lambda_s \notin \mathbb{Z}$ and 
$\sum_{i=1}^{s} \lambda_i =0$, 
we have 
\begin{enumerate}
	\item 
$H^k(M(\A), \tilde{\mathcal{L}}_{\lambda}) = 0$ 
for $k \not= n, n-1$, 
	\item 
$H^{n}(M(\A), \tilde{\mathcal{L}}_{\lambda})  
 \simeq H^{n-1}(M(\A), \tilde{\mathcal{L}}_{\lambda})$ 
and $\dim{H^{n}} = \dim{H^{n-1}} = \binom{s-2}{n-1}$.  
	\item 
$\{ e_1 \we e_{i_1} \we \cdots \we e_{i_{n-1}}  
: 1 < i_1 < \cdots < i_{n-1} < s \}$ is a basis of 
$H^{n}$. 
	\item 
$\{ \der{[ e_1 : e_{i_1} : \cdots : e_{i_{n-1}} ]} 
: 1 < i_1 < \cdots < i_{n-1} < s \}$ is a basis of 
$H^{n-1}$, 
\end{enumerate}
where $e_j = \dlogg{\alpha_j}$ and 
$\tilde{\mathcal{L}}_{\lambda}$ is the rank one local system on $M(\A)$ 
whose monodromy around the hyperplane $H_j$ is 
$\exp{(-2 \pi \sqrt{-1} \lambda_j)}$.  
\end{lem} 
Note that 
$H^k(M(\A), \tilde{\mathcal{L}}_{\lambda})$ is isomorphic to 
the twisted de Rham cohomology defined by the one form  
$e_{\lambda} = \sum_{j=1}^s \lambda_j e_j$ (see \cite{OT2}) 
and  that 
$\omega_{\lambda} \vert_U = e_{\lambda}$ 
and 
$\eta{[i_1, \ldots, i_{n}]} \vert_U = 
\der{[e_{i_1} : \cdots : e_{i_{n}}]}$.  
Now suppose that 
there exists a global ($n-2$)-form $\alpha$ such that 
$\eta{[i_1, \ldots, i_{n}]} = \nabla_{\lambda}(\alpha)$. 
Then restricting it to $U$, 
we have 
$\der{[ e_{i_1} : \cdots : e_{i_{n}}]} 
= \tilde{\nabla}_{\lambda}(\alpha \vert_U)$ 
where $\tilde{\nabla}_{\lambda} = d + e_{\lambda} \we$. 
However, by the above Lemma, this is a contradiction. 
Therefore   
$\eta{[i_1, \ldots, i_{n}]}$ defines a non-vanishing class of degree $n-1$ 
for $1 \leq i_1 < \cdots < i_{n} \leq s$.  
In a similar fashion, 
we obtain 
$\{ \eta{[1, i_1, \ldots, i_{n-1}]}  
: 1 < i_1 < \cdots < i_{n-1} < s \}$ is independent 
in $H^{n-1}(M, \mathcal{L}_{\lambda})$.   

Assume $s<m$ and fix $m$. 
Take $\eta{[m, i_1, \ldots, i_{n}]}$ for $1 \leq i_1 < \cdots < i_{n} \leq s$.  
It is easy to see that  
$\eta{[m, i_1, \ldots, i_{n}]} \vert_U = e_{i_1} \we \cdots \we e_{i_n}$. 
Therefore, 
$\eta{[m, i_1, \ldots, i_{n}]}$ defines a non-vanishing class of degree $n$.  
By the same way,  using the above Lemma, 
we have 
$\{ \eta{[m, 1,  i_1, \ldots, i_{n-1}]}  
: 1 < i_1 < \cdots < i_{n-1} < s \}$ is independent 
in $H^{n}(M, \mathcal{L}_{\lambda})$.   
This completes the proof.   
\end{proof}

If a weight $\lambda$ is trivial then 
$H^{k}(M, \mathcal{L}_{\lambda})$ is isomorphic to 
the usual de Rham cohomology $H^{k}(M)$ on $M$.  

\begin{cor} 
Under the assumption of Theorem \ref{MainTheorem}, 
we have 
$$   \dim{H^{k}(M)} \geq  \binom{s-1}{k}   
\quad \text{ for } 1 \leq k \leq n-1.   $$ 
Moreover, if $s < m$ then we have 
$$  \dim{H^{n}(M)} \geq  \binom{s-1}{n-1}  \  \text{ and }   \  
    \dim{H^{k}(M)} \geq  \binom{s}{k} \ \text{ for } 1 \leq k \leq n-1.  $$ 
\end{cor} 

\begin{proof} 
In the proof of Theorem \ref{MainTheorem},  
since $\A$ is the generic arrangement of $s$ hyperplanes in $\mathbb{C}^n$ 
and $\dc{\A}$ is in general position, 
the following is known (\cite{OT}). 
\begin{enumerate}
	\item 
$\dim{H^k(M(\A))} = \binom{s}{k}$ for $1 \leq k \leq n$ and 
$\dim{H^k(M(d \A))} = \binom{s-1}{k}$ for $1 \leq k \leq n-1$,  
	\item 
$H^k(M(\A)) \simeq H^k(M(d \A)) \oplus  H^{k-1}(M(d \A))$ for $1 \leq k \leq n-1$,    
 and $H^n(M(\A)) \simeq H^{n-1}(M(d \A))$, 
	\item 
$\{ e_1 \we e_{i_1} \we \cdots \we e_{i_{n-1}}  
: 1 < i_1 < \cdots < i_{n-1} \leq s \}$ is a basis of $H^{n}(M(\A))$. 
	\item 
$\{ \der{[ e_1 : e_{i_1} : \cdots : e_{i_{k}} ]} 
: 1 < i_1 < \cdots < i_{k} \leq s \} 
\cup \{ e_1 \we e_{i_1} \we \cdots \we e_{i_{k-1}}  
: 1 < i_1 < \cdots < i_{k-1} \leq s \}$ 
is a basis of 
$H^k(M(\A))$ for $1 \leq k \leq n-1$. 
\end{enumerate} 
Since $\eta{[i_1, \ldots, i_{k}]}$ is $d$-closed, 
the same argument leads this corollary. 
\end{proof}

\section{Generalization}\label{generalization} 

Let $\A$ be the arrangement of (affine or projective) hyperplanes. 
The \textit{intersection set} $L(\A)$ of $\A$ is 
the set of nonempty intersections of elements of $\A$. 
For $X \in L(\A)$, define a central arrangement  
$\A_{X} = \{ H \in \A \ \vert \ X \subset H \}$.  
Let $\C$ be a central arrangement 
with center $\bigcap_{H \in \C} H \not= \emptyset$. 
We call $\C$ \textit{decomposable} 
if there exist nonempty subarrangements $\C_1$ and $\C_2$ 
so that $\C=\C_1 \cup \C_2$ is a disjoint union and 
after a linear coordinate change 
the defining polynomials for $\C_1$ and $\C_2$   
have no common variables. 
Define 
$\den(\A) = \{ X \in L(\A) : \A_{X} \text{ is not decomposable.} \}$.  
For a complex weight $\lambda$ of $\A$ and $X \in L(\A)$,  
denote $\lambda_{X} = \sum_{H \in \A_X} \lambda_{H}$. 
The construction of a basis of the twisted cohomology for arrangements 
given by \cite{FT} (cf. \cite{OT2}), 
by the same way of proof of Theorem \ref{MainTheorem},  
induces the following.  

\begin{thm} 
Let $1 < n < s \leq m$. 
Let $D_1, \ldots, D_m$ be effective divisors on $\mathbb{P}^n$ 
with same degree   
and let $M$ be the complement of their supports.  
Suppose (A1) and (A2). 
Let $\A$ be the arrangement of hyperplanes in 
the dual projective space 
$\Lambda^{*} = (\mathbb{P}^{n-1})^* \simeq \mathbb{P}^{n-1}$ 
defined by $D_1, \ldots, D_s$.  
Let $\lambda$ be a non-trivial weight 
such that $\sum_{i=1}^s \lambda_i =0$ 
and $\lambda_i =0$ for $i = s+1, \ldots, m$.  
If $\lambda_{X} \not\in \mathbb{Z}_{\geq 0}$ 
for every $X \in \den(\A)$, 
then we have 
$$   \dim{H^{n-1}(M, \mathcal{L}_{\lambda})} \geq  \beta,    $$ 
and moreover, if $s < m$ then we have 
$$   \dim{H^{n}(M, \mathcal{L}_{\lambda})} \geq  \beta,   $$ 
where $\beta$ is the Euler characteristic $\chi(M(\A))$ of 
$M(\A) = \mathbb{P}^{n-1} \setminus \cup_{H \in \A} H$. 
\end{thm}
\begin{rem} 
Note that $\beta$ is known as 
the beta invariant of the underlying matroid of $\A$. 
If $\A$ is defined over real then $\beta$ is 
the number of bounded chambers in 
$\mathbb{C}^{n-1} = \mathbb{P}^{n-1} \setminus H$ for fixed $H \in \A$ 
(see \cite{STV, OT2}).   
\end{rem} 
\begin{rem} 
Note that 
$\lambda_H \not\in \mathbb{Z}_{\geq 0}$ for all hyperplanes in $\A$.  
We can generalize this theorem 
in the case that there is $H \in \A$ with $\lambda_H =0$, 
by using \cite{Ka2}. 
\end{rem}

\section{Affine case}\label{affine} 

Let $V_1^a, \ldots, V_m^a$ be hypersurfaces 
in the complex affine space $\mathbb{C}^n$ 
with coordinates $u= (u_1, \ldots, u_n)$. 
Denote $V^a= \cup_{i=1}^m V_i^a$ and $M = \mathbb{C}^n \setminus V^a$. 
Assume that a polynomial $f_j(u)$  of degree $d_i$ defines $V_j^a$. 
Let $\lambda = (\lambda_1, \ldots, \lambda_m)$ be a weight and 
let $\omega_{\lambda}^a = \sum_{i=1}^{m} \lambda_i \dlogg{f_i}$. 
Then we obtain the twisted de Rham complex  
$(\Omega(* V^a), \nabla_{\lambda}^a)$ 
where  $\Omega(* V^a)$ is the space of rational forms with poles along $V^a$ 
and  $\nabla_{\lambda}^a = d + \omega_{\lambda}^a \we$. 
The Grothendiek-Deligue comparison theorem (\cite{De}) asserts  
$$   H^{k}(M, \mathcal{L}_{\lambda}^a)  \simeq 
     H^k(\Omega(* V^a), \nabla_{\lambda}^a),   $$ 
where $\mathcal{L}_{\lambda}^a$ is the rank one local system defined by 
the flat connection $\nabla_{\lambda}^a$ (cf. \cite{KN, Ki}).  

Let $\mathbb{P}^n$ be the complex projective space 
with the infinite hyperplane $H_{\infty}$, 
which is a compactification of $\mathbb{C}^n$. 
Let $[x_0: \ldots: x_n]$ denote homogeneous coordinates 
with $H_{\infty} = \{ x_0 =0 \}$. 
Define the homogeneous polynomial 
$F_j (x) =  x_0^d f_j( x_1/x_0, \ldots, x_n/x_0)$ 
of degree $d = \max{(d_1, \ldots, d_m)}$.  
Then $F_j (x)$ determines the divisor $D_j$ of degree $d$.   
If $d = d_1 = \cdots = d_m$ then 
the support of $\sum_{i=1}^m D_i$ does not contain $H_{\infty}$,  
otherwise it contains $H_{\infty}$.  
Note that the weight of $H_{\infty}$ is given by 
$- \sum_{i=1}^m \lambda_i d_j$.   
If $d = d_1 = \cdots = d_m$ and $\sum_{i=1}^m \lambda_i=0$ then 
the weight of $H_{\infty}$ is zero.  

\begin{cor} \label{Cor:affine}
Under the assumption in Theorem  \ref{MainTheorem},  
if a point $P$ in (A2) is not on $H_{\infty}$, then we have 
$$   \dim{H^{n-1}(M^a, \mathcal{L}_{\lambda}^a)} \geq  \binom{s-2}{n-1}  $$ 
for a non-trivial weight $\lambda$ with $\sum_{i=1}^s \lambda_i=0$ and 
$\lambda_i =0$ for $i = s+1, \ldots, m$. 
\end{cor}

\section{Examples} \label{examples} 

\subsection{$n=2$ and $s=3$} 
Let $D_1$, $D_2$ and $D_3$ be irreducible prime divisors on $\mathbb{P}^2$ 
defined by 
$F_1 = x_0^2 + x_1^2 - 2 x_2^2$, 
$F_2 = x_0^2 + 2 x_1^2 -3 x_2^2$ and  
$F_3 = 2x_0^2 + x_1^2 -3 x_2^2$, respectively. 
They are three generic elements of the pencil of conic curves 
$F_{[a:b]} = a (x_0^2 - x_2^2) +  b (x_1^2 - x_2^2)$ 
and 
they intersect transversally each other at each four intersection points. 
Let $M = \mathbb{P}^2 \setminus \cup_{i=1}^3 D_i$ and 
$M^a = \mathbb{C}^2 \setminus \cup_{i=1}^3 (D_i \cap \mathbb{C}^2)$ 
where $\mathbb{C}^2 = \mathbb{P}^2 \setminus \{ x_2=0 \}$.  
We have 
$H^1(M, \mathcal{L}_{\lambda})  \not= 0$ and 
$H^1(M^a, \mathcal{L}_{\lambda}) \not= 0$ 
for a non-trivial weight $\lambda$ with $\lambda_1 + \lambda_2 + \lambda_3=0$. 
In the following case, we can get it also. 
\begin{enumerate}
	\item 
$F_1 = x_0^2 + x_1^2 - 2 x_2^2$, 
$F_2 = x_0^2 + 2 x_1^2 -3 x_2^2$ and  
$F_3' = x_0^2 -x_1^2$ (two conics and a set of 2-lines).  
	\item 
$F_1 = x_0^2 + x_1^2 - 2 x_2^2$, 
$F_2' =  x_2^2 - x_0^2$ and  
$F_3' = x_0^2 -x_1^2$ (one conics and two sets of 2-lines).  
	\item 
$F_1' = x_1^2 - x_2^2$, 
$F_2' = x_2^2 - x_0^2$ and  
$F_3' = x_0^2 -x_1^2$ (three sets of 2-lines).  
\end{enumerate}
In the last case, it is  
an arrangement of 6 lines in the Ceva Theorem.  
In a similar fashion, we get the following. 
In the degree three case, 
we get an arrangement of 9 lines in the Pappus Theorem (\cite{Fa}).  
In the degree four case, 
there are two different arrangements of 12 lines in  
the Kirkman Theorem and the Steiner Theorem (\cite{Ka3}). 
Those arrangements are 3-nets, 
whose combinatorial structures are  
matroids associated to Latin squares (\cite{LY, Yu2, Ka3}).  
In the other hand, the $B_3$-arrangement is 
an example of the case that $D_i$'s are not prime. 
Let $D_1$, $D_2$ and $D_3$ be divisors  
defined by 
$x_2^2(x_0^2 - x_1^2)$, 
$x_1^2(x_0^2 - x_2^2)$ and   
$x_0^2 (x_1^2 - x_2^2)$, respectively. 
Their divisors can be written by  
$D_1 = 2 H_1 + H_2 + H_3$, 
$D_2 = 2 H_4 + H_5 + H_6$ and 
$D_3 = 2 H_7 + H_8 + H_9$  
where $H_i$'s are hyperplanes.   
The arrangement $\A = \{ H_1, \ldots, H_9 \}$ 
is called the $B_3$-arrangement. 
Note that a weight $\lambda$ induces the weight 
$(2 \lambda_1, \lambda_1, \lambda_1, 2 \lambda_2, \lambda_2, \lambda_2, 
2 \lambda_3, \lambda_3, \lambda_3 )$ of $\A$ (cf. \cite{Fa, Ka3}). 

\subsection{$n=2$ and $s \geq 3$} 
We consider the pencil of cubic curves 
$F_{[a:b]} = a (x_0^3 + x_1^3 + x_2^3) + 3 b x_0 x_1 x_2$.  
A generic element given by $a \not= 0$ and $b^3 \not= -1$ 
is non-singular.  
For $s$ generic elements $D_1, \ldots, D_s$, 
we have  $\dim{H^1(M, \mathcal{L}_{\lambda})} \geq  s-2$.   
Define non-generic elements 
$F_1 = x_0 x_1 x_2$, 
$F_2 = x_0^3 + x_1^3 + x_2^3 -3 x_0 x_1 x_2$, 
$F_3 = x_0^3 + x_1^3 + x_2^3 -3 \xi x_0 x_1 x_2$ and 
$F_4 = x_0^3 + x_1^3 + x_2^3 -3 \xi^2 x_0 x_1 x_2$  
where  $\xi = e^{2 \pi \sqrt{-1} /3}$.  
They are four sets of 3-lines and 
the set $\A$ of all lines is called 
the Hessian configuration, 
which is the arrangement of 12 lines 
passing through the nine inflection points of a nonsingular cubic.   
In this case, we know $\dim{H^1(M, \mathcal{L}_{\lambda})} = 2$ 
for a non-trivial weight $\lambda$ with $\sum_{i=1}^{4} \lambda_i =0$ 
(\cite{Li}).    

\subsection{$n=3$} 
Let 
$F_1 = x_0 (      x_1   + x_2  + x_3)$, 
$F_2 = x_1 (-x_0        + x_2  - x_3)$, 
$F_3 = x_2 (-x_0 - x_1         + x_3)$ and 
$F_4 = x_3 (-x_0 + x_1  - x_2       )$.  
Then we can check divisors defined them satisfy 
the conditions in Theorem \ref{MainTheorem} 
and then  $H^{2}(M, \mathcal{L}_{\lambda}) \not= 0$  
for a weight $\lambda$ with  $\sum_{i=1}^{4} \lambda_i =0$. 
They yield the arrangement of 8 planes, 
whose underlying matroid is of type $L_8$ (\cite{Ka3}).

\subsection{higher dimensional case} 
Let $F_i = x_{i -1}^d - x_{i}^d$ for $i = 1, \ldots, n$ 
and $F_0 = x_n^d - x_0^d$. 
This support determines the arrangement $\A$ 
of $(n+1)d$ hyperplanes in $\mathbb{P}^n$. 
Note that this is a projective closure of a subarrangement of 
the monomial arrangement $\A_{d,d,n+1}$ (see \cite{OT, CS}). 
Since $\sum_{k=0}^{n} F_k =0$, we have 
$H^{n-1}(M, \mathcal{L}_{\lambda}) \not= 0$  
for a weight $\lambda$ with $\sum_{i=0}^{n} \lambda_i =0$.  
The $n=2$ case was found in  \cite{CS}.  
Note that the underlying matroid of $\A$ is  
a degeneration of the matroid associated to 
the Latin $n$-dimensional hypercube  
given by the addition table for $(\mathbb{Z}_d)^n$ (see \cite{Ka3}).

 
\end{document}